# EDGE-ENHANCING RECONSTRUCTION ALGORITHM FOR THREE-DIMENSIONAL ELECTRICAL IMPEDANCE TOMOGRAPHY

L. HARHANEN[†], N. HYVÖNEN[‡], H. MAJANDER[‡], AND S. STABOULIS[‡]

**Abstract.** Electrical impedance tomography is an imaging modality for extracting information on the conductivity distribution inside a physical body from boundary measurements of current and voltage. In many practical applications, it is *a priori* known that the conductivity consists of embedded inhomogeneities in an approximately constant background. This work introduces an iterative reconstruction algorithm that aims at finding the *maximum a posteriori* estimate for the conductivity assuming an edge-preferring prior. The method is based on applying (a single step of) priorconditioned lagged diffusivity iteration to sequential linearizations of the forward model. The algorithm is capable of producing reconstructions on dense unstructured three-dimensional finite element meshes and with a high number of measurement electrodes. The functionality of the proposed technique is demonstrated with both simulated and experimental data in the framework of the complete electrode model, which is the most accurate model for practical impedance tomography.

**1. Introduction.** The aim of *electrical impedance tomography* is to reconstruct the internal conductivity distribution of a physical body based on boundary measurements of current and voltage. This constitutes a nonlinear and severely illposed inverse problem. EIT can be used in medical imaging, process tomography and non-destructive testing of materials. Consult the review articles [2, 10, 33] for more information on EIT, the associated mathematical theory and the related reconstruction algorithms. In this work, we consider EIT under the prior assumption that the to-be-reconstructed conductivity consists of well defined inclusions in an approximately homogeneous background, which is a setting encountered in many practical applications: Consider, e.g., the localization of air bubbles or manufacturing defects in a piece of building material. We work exclusively with the *complete electrode model* (CEM), which is the most accurate model for real-world EIT [11, 32]; in particular, the introduced algorithm also estimates the contact resistances that are an unavoidable nuisance of practical EIT.

We tackle the reconstruction problem of EIT within the Bayesian paradigm and incorporate the prior information on the structure of the imaged object by introducing an edge-preferring prior density for the (discretized) conductivity. The prior for the contact resistances is chosen to be uninformative since their estimation from EIT measurements is not an illposed problem. Assuming an additive Gaussian measurement noise model, the computation of the *maximum a posteriori* (MAP) estimate, i.e., the maximizer of the posterior density, corresponds to finding a minimizer for a Tikhonov-type functional that exhibits nonquadratic behavior in both the discrepancy and the penalty term. In the terminology of regularization theory, this corresponds to employing Tikhonov regularization with a suitable nonlinear penalty term, such as *total*

[†]Aalto University, Department of Engineering Design and Production, P.O. Box 11000, FI-00076 Aalto, Finland (lauri.harhanen@aalto.fi).

[‡]Aalto University, Department of Mathematics and Systems Analysis, P.O. Box 11100, FI-00076 Aalto, Finland (nuutti.hyvonen@aalto.fi, helle.majander@aalto.fi, stratos.staboulis@aalto.fi). The work of NH, HM and SS was supported by the Academy of Finland (decision 267789).

*variation* (TV) [31] or *Perona–Malik* (PM) [30]. For previous works employing edge-enhancing or sparsity-promoting regularization/prior in two- and three-dimensional EIT, see, e.g., [3, 12, 18, 24] and [26], respectively.

The goal of this work is to introduce and numerically test a novel algorithm for computing EIT reconstructions by minimizing the aforementioned Tikhonov functional, stressing reconstructions on dense and unstructured *three-dimensional* finite element meshes; see [1] for a closely related technique for linear inverse problems. In its basic form, the method is composed of two nested iterations: The outer loop corresponds to sequential linearizations of the CEM, while the inner loop minimizes the resulting Tikhonov functionals, which correspond to (large) illconditioned matrix equations with nonquadratic penalty terms for the conductivity, by means of the *lagged diffusivity* iteration [37]. However, the final version of our algorithm is based on just a single lagged diffusivity step for each linearization of the forward model, which is computationally efficient and does not seem to impair the reconstruction quality. The use of sequential linearizations is motivated by the theoretical observation that — under the so-called *continuum model* — the inclusion shapes are not affected if the nonlinear inverse problem of EIT is replaced by its linearized version [21].

In our framework, a lagged diffusivity step corresponds to evaluating the coefficient matrix appearing in the necessary condition for the minimizer of the considered Tikhonov functional (with a nonquadratic penalty term) at the previous iterate when computing the next one (cf. (4.5)). The resulting linear equation can be interpreted as the one defining the MAP estimate for the linearized inverse problem of the CEM with an improper Gaussian prior: The inverse covariance matrix $H$ for the conductivity originates from the gradient of the nonlinear penalty term at the previous step whereas the prior for the contact resistances is still uninformative. In particular, the equations for the conductivity and the contact resistances can be decoupled, which enables applying *priorconditioning* [5, 6, 7, 9] solely to the conductivity.

We precondition the aforementioned matrix equation for the conductivity by $H$; since the preconditioner is the ('linearized') inverse prior covariance matrix, such a procedure corresponds to whitening the prior and was dubbed priorconditioning in [5, 6, 7, 9]. Employing LSQR [29, 28] on the priorconditioned equation leads to searching for the new lagged diffusivity iterate for the conductivity from Krylov subspaces that have the prior information related to $H$ as a built-in property. This accelerates the LSQR iteration considerably (cf. [1]). The computation of a symmetric factorization for the preconditioner $H$ is avoided by resorting to the LSQR version introduced in [1], which makes the algorithm particularly attractive for unstructured finite element meshes that are needed for the CEM due to, e.g., the singularities at the electrode edges (cf. [15, 32]).

Our numerical experiments consider exclusively three-dimensional settings. It is demonstrated that the proposed algorithm produces good quality reconstructions on dense FEM meshes in computational times of just a few minutes on a standard desktop computer. A couple of different geometrical settings are considered, and reconstructions are computed from simulated noisy electrode measurements as well as from experimental water tank data.

This text is organized as follows. Section 2 recalls the CEM and Section 3 introduces the discrete Bayesian setting for the inverse problem of EIT. In Section 4, we describe the basic ideas behind the reconstruction algorithm, while Section 5 considers implementation aspects. The numerical experiments and concluding remarks are presented in Sections 6 and 7, respectively.

**2. Complete electrode model.** Let $\Omega \subset \mathbb{R}^n$, $n = 2$ or $3$, be a bounded domain with a smooth enough boundary and $\sigma \in L^\infty(\Omega)$ the corresponding conductivity such that $\sigma \geq c > 0$ almost everywhere in $\Omega$. The object boundary $\partial\Omega$ is assumed to be partially covered by $M \in \mathbb{N} \setminus \{1\}$ connected electrodes $\{E_m\}_{m=1}^M$ that do not overlap and are identified with open subsets of $\partial\Omega$. We denote $E = \cup_m E_m$. The electrode net current and voltage patterns are represented by $I, U \in \mathbb{R}_\diamond^M$, respectively, with $I_m, U_m \in \mathbb{R}$ corresponding to the measurements on the $m$th electrode. Here, $\mathbb{R}_\diamond^M \subset \mathbb{R}^M$ is the subspace of mean-free vectors; the current patterns belong to $\mathbb{R}_\diamond^M$ due to conservation of charge, the electrode voltages via choosing the ground level of potential appropriately. The contact resistances (cf. [11]) that characterize the thin and highly resistive layers at the electrode-object interfaces are modelled by $z \in \mathbb{R}^M$ with $z_m \geq c > 0$, $m = 1, \ldots, M$.

Let $u$ denote the electromagnetic potential inside $\Omega$. According to the CEM [11, 32], the pair $(u, U) \in H^1(\Omega) \oplus \mathbb{R}_\diamond^M$ is the unique solution of the elliptic boundary value problem

$$\begin{aligned}
\nabla \cdot (\sigma \nabla u) &= 0 & &\text{in } \Omega, \\
\frac{\partial u}{\partial \nu} &= 0 & &\text{on } \partial\Omega \setminus \overline{E}, \\
u + z_m \sigma \frac{\partial u}{\partial \nu} &= U_m & &\text{on } E_m, \quad m = 1, \ldots, M, \\
\int_{E_m} \sigma \frac{\partial u}{\partial \nu} \, dS &= I_m, & &m = 1, \ldots, M,
\end{aligned} \qquad (2.1)$$

for a given electrode current pattern $I \in \mathbb{R}_\diamond^M$ and with $\nu : \partial\Omega \to \mathbb{R}^n$ denoting the exterior unit normal of $\partial\Omega$. For a physical justification of (2.1), we refer to [11].

We denote by $U(\sigma, z; I) \in \mathbb{R}_\diamond^M$ the functional dependence of the electrode potentials on the conductivity, the contact resistances and the applied electrode currents. This notation is extended to the case of $M-1$ linearly independent current patterns, $I^1, \ldots, I^{M-1} \in \mathbb{R}_\diamond^M$ by writing

$$\mathcal{U}(\sigma, z; I^1, \ldots, I^{M-1}) = \left[U(\sigma, z; I^1)^{\mathrm{T}}, \ldots, U(\sigma, z; I^{M-1})^{\mathrm{T}}\right]^{\mathrm{T}} \in \mathbb{R}^{M(M-1)},$$

where the dependence on $I^1, \ldots, I^{M-1}$ is usually omitted for the sake of clarity. The inverse problem of EIT consists of reconstructing (useful information about) the conductivity $\sigma$ — and the contact resistances $z$ — from (a noisy version) of the electrode potential measurements $\mathcal{U}(\sigma, z; I^1, \ldots, I^{M-1})$.

**3. Bayesian framework and edge-preserving priors.** In the rest of this work, we consider a discretized version of (2.1). To be more precise, the conductivity $\sigma$ is modelled as

$$\sigma = \sum_{j=1}^N \sigma_j \varphi_j, \qquad (3.1)$$

where $\sigma_j \in \mathbb{R}_+$, $j = 1, \ldots, N$, and $\varphi_j \in H^1(\Omega)$, $j = 1, \ldots, N$, are the piecewise linear basis functions corresponding to a finite element mesh of $\Omega$. We denote by $\sigma$ both the vector $\sigma \in \mathbb{R}_+^N$ and the corresponding conductivity defined by (3.1); the exact meaning should be clear from the context. Moreover, by the solution $(u, U)$ of (2.1) we mean the approximate FEM solution in the subspace

$$\operatorname{span}\{\varphi_j\}_{j=1}^N \oplus \mathbb{R}_\diamond^M \subset H^1(\Omega) \oplus \mathbb{R}_\diamond^M;$$

see [34, 35] for the implementation details.

The potentials measured at the electrodes are modelled as

$$\mathcal{V} = \mathcal{U}(\sigma, z) + \eta, \tag{3.2}$$

where $\eta \in \mathbb{R}^{M(M-1)}$ is a realization of a Gaussian random variable with zero mean and a known, symmetric and positive definite covariance matrix $\Gamma \in \mathbb{R}^{M(M-1) \times M(M-1)}$. In consequence, the *likelihood function*, i.e., the probability density of the measurement $\mathcal{V}$ given the parameters $\sigma$ and $z$, is

$$p(\mathcal{V} \,|\, \sigma, z) \,\propto\, \exp\Big( -\frac{1}{2} \big(\mathcal{V} - \mathcal{U}(\sigma, z)\big)^{\mathrm{T}} \Gamma^{-1} \big(\mathcal{V} - \mathcal{U}(\sigma, z)\big) \Big).$$

Here and in what follows, we denote by $p$ a generic probability density whose exact meaning should be clear from the context.

In Bayesian inversion the unknown parameters $\sigma$ and $z$ are given prior probability densities that reflect the available information before the measurements. We assume that the (discretized) conductivity is *a priori* distributed according to the density

$$p(\sigma) \propto \exp\big( -\alpha R(\sigma) \big), \tag{3.3}$$

where $\alpha > 0$ is a free parameter and $R$ is of the form

$$R(\sigma) = \int_\Omega r\big(|\nabla \sigma(x)|\big) \, dx, \tag{3.4}$$

with $r : \Omega \to \mathbb{R}_+$ being a suitable, continuously differentiable, monotonically increasing function that prefers edges over slow changes in the conductivity $\sigma$. (In Section 5, we will complement the prior density (3.3) with assumptions on the boundary behavior of $\sigma$ and also get rid of $\alpha$.) In this work, we consider the total variation and Perona–Malik priors,

$$r(t) = \sqrt{T^2 + t^2} \approx |t| \qquad \text{and} \qquad r(t) = \frac{1}{2} T^2 \log\big(1 + (t/T)^2\big), \tag{3.5}$$

where $T > 0$ is a small parameter that ensures $r$ is differentiable (TV) [31] or controls a threshold for detectable edges (PM) [30]. On the other hand, the contact resistances are given an *uninformative* prior, that is, all realizations of $z$ in $\mathbb{R}^M$ are considered equally probable *a priori*. The reason for such a generic choice is that, according to our experience, the determination of the contact resistances from EIT measurements is not an illposed problem, and thus no prior knowledge (or regularization) is needed.

By the Bayes' formula, the posterior density, i.e., the density for the pair $(\sigma, z) \in \mathbb{R}^N \times \mathbb{R}^M$ given the measurement $\mathcal{V}$, is

$$\begin{aligned} p(\sigma, z \,|\, \mathcal{V}) &\propto p(\mathcal{V} \,|\, \sigma, z) \, p(\sigma) \\ &\propto \exp\Big( -\frac{1}{2} \big(\mathcal{V} - \mathcal{U}(\sigma, z)\big)^{\mathrm{T}} \Gamma^{-1} \big(\mathcal{V} - \mathcal{U}(\sigma, z)\big) - \alpha R(\sigma) \Big), \end{aligned}$$

where the constants of proportionality do not depend on $\sigma$ and $z$. In this work, we propose an efficient algorithm for (approximately) determining the (possibly nonunique) MAP estimate corresponding to this posterior, that is, the minimizer of the Tikhonov functional

$$\Phi(\sigma, z) := \frac{1}{2} \big(\mathcal{V} - \mathcal{U}(\sigma, z)\big)^{\mathrm{T}} \Gamma^{-1} \big(\mathcal{V} - \mathcal{U}(\sigma, z)\big) + \alpha R(\sigma). \tag{3.6}$$

Notice that $\Phi$ exhibits nonquadratic behavior in the term originating from the likelihood function as well as in the one penalizing for the unwanted behavior of $\sigma$. In the following, we denote $\kappa = [\sigma^\mathrm{T}, z^\mathrm{T}]^\mathrm{T} \in \mathbb{R}^{N+M}$ and slightly abuse the notation by sometimes writing $\Phi(\kappa)$ and $\mathcal{U}(\kappa)$ instead of $\Phi(\sigma, z)$ and $\mathcal{U}(\sigma, z)$, respectively.

REMARK 3.1. *One would end up minimizing a functional of the form* (3.6) *by just applying Tikhonov regularization with an edge-preferring, nonquadratic penalty term to the reconstruction problem EIT within the CEM. In fact, such an interpretation is arguably more natural for many considerations in the following sections; for example, the employed priorconditioning can be considered as the transformation of a (quadratic) Tikhonov functional into the standard form, which is a well-established idea* [16, 20, 22]. *However, since priorconditioning, i.e., preconditioning with the inverse prior covariance matrix, is motivated by Bayesian ideas, the chosen statistical framework is well-grounded.*

**4. Sequential linearization and lagged diffusivity iterations.** The density (3.3) reflects the prior information that the target conductivity consists of well defined inclusions in approximately homogeneous background. On the other hand, it is known that the reconstruction of inclusion supports is, loosely speaking, unaffected if the nonlinear inverse problem of EIT is replaced by its linearized version [21]. These observations give the motivation for our preliminary plan of minimizing $\Phi$ by an iterative two-level algorithm: Start from an initial guess $\kappa^{(0)} = [(\sigma^{(0)})^\mathrm{T}, (z^{(0)})^\mathrm{T}]^\mathrm{T}$ and set $j = 0$.

1. Replace the nonlinear forward model $\mathcal{U}(\kappa)$ in (3.6) by its linearization around the current iterate $\kappa^{(j)} = [(\sigma^{(j)})^\mathrm{T}, (z^{(j)})^\mathrm{T}]$ and write the necessary condition for the minimizer.
2. Minimize the resulting Tikhonov functional by the *lagged diffusivity iteration* [37]. Denote the minimizer as $\kappa^{(j+1)}$, set $j \leftarrow j+1$ and return to step 1, unless the chosen stopping criterion is satisfied.

We discuss the details of these two steps in the following subsections. It should be mentioned, however, that the final implementation of our reconstruction algorithm is a simplified version of the above scheme in order to optimize the computational efficiency; see Section 5.

**4.1. Linearization of the forward model.** The linearization of $\mathcal{U}(\kappa)$ around $\kappa^{(j)}$ in (3.6) leads to a standard Tikhonov functional with a nonquadratic penalty term,

$$\Phi^{(j)}(\kappa) := \frac{1}{2}\big(y^{(j)} - J^{(j)}\kappa\big)^\mathrm{T} \Gamma^{-1} \big(y^{(j)} - J^{(j)}\kappa\big) + \alpha\, R(\sigma), \qquad (4.1)$$

where $\kappa = [\sigma^\mathrm{T}, z^\mathrm{T}]^\mathrm{T}$, the matrix $J^{(j)} \in \mathbb{R}^{M(M-1) \times (N+M)}$ is the Jacobian of the map $\kappa \mapsto \mathcal{U}(\kappa)$ evaluated at $\kappa^{(j)}$ and

$$y^{(j)} = \mathcal{V} - \mathcal{U}(\kappa^{(j)}) + J^{(j)}\kappa^{(j)} \ \in \mathbb{R}^{M(M-1)}.$$

See, e.g., [24, 36] for instructions on forming the Jacobian $J^{(j)}$. The necessary condition for a (local) minimizer of (4.1) is $\nabla \Phi^{(j)} = 0$, which yields the equation

$$(J^{(j)})^\mathrm{T} \Gamma^{-1} J^{(j)} \kappa + \alpha \begin{bmatrix} (\nabla R)(\sigma) \\ 0 \end{bmatrix} = (J^{(j)})^\mathrm{T} \Gamma^{-1} y^{(j)}, \qquad (4.2)$$

where $0 \in \mathbb{R}^M$ denotes a zero vector of the appropriate length.

It is straightforward to check that the partial derivatives of the penalty term $R : \mathbb{R}^N \to \mathbb{R}_+$ are given by the formula

$$\frac{\partial}{\partial \sigma_k} R(\sigma) = \int_\Omega \frac{r'(|\nabla \sigma(x)|)}{|\nabla \sigma(x)|} \nabla \sigma(x) \cdot \nabla \varphi_k(x) \, dx$$

$$= \sum_{l=1}^N \left( \int_\Omega \frac{r'(|\nabla \sigma(x)|)}{|\nabla \sigma(x)|} \nabla \varphi_k(x) \cdot \nabla \varphi_l(x) \, dx \right) \sigma_l, \qquad k = 1, \ldots, N,$$

where $\varphi_j$ is the $j$th piecewise linear FEM basis function and $\sigma : \Omega \to \mathbb{R}_+$ is defined by (3.1). Introducing the $\sigma$-dependent matrix

$$H_{k,l}(\sigma) := \int_\Omega \frac{r'(|\nabla \sigma(x)|)}{|\nabla \sigma(x)|} \nabla \varphi_k(x) \cdot \nabla \varphi_l(x) \, dx, \qquad k, l = 1, \ldots, N, \qquad (4.3)$$

the gradient of $R$, needed in (4.2), allows the representation

$$(\nabla R)(\sigma) = H(\sigma) \sigma.$$

Denoting

$$c_\sigma : x \mapsto \frac{r'(|\nabla \sigma(x)|)}{|\nabla \sigma(x)|}, \quad \Omega \to \mathbb{R}_+,$$

$H(\sigma) \in \mathbb{R}^{N \times N}$ can be interpreted as the FEM system matrix for the elliptic partial differential operator

$$-\nabla \cdot c_\sigma \nabla \qquad (4.4)$$

accompanied with the natural boundary condition on $\partial \Omega$. Obviously $H(\sigma)$ is not invertible since neither the operator (4.4) itself or the natural boundary condition, i.e., a Neumann condition with a to-be-specified right-hand side, sees the 'ground level of potential'. We fix this problem by accompanying (4.4) with a suitable Dirichlet boundary condition on a subset of $\partial \Omega$. The details can be found in Section 5; in the rest of this section, we just assume that the inverse $H^{-1}$ exists.

Let us write $J^{(j)} = [J_1^{(j)}, J_2^{(j)}]$, where $J_1^{(j)} \in \mathbb{R}^{M(M-1) \times N}$ contains the derivatives with respect to $\sigma \in \mathbb{R}^N$ and $J_2^{(j)} \in \mathbb{R}^{M(M-1) \times M}$ those with respect to $z \in \mathbb{R}^M$. According to the above considerations, the necessary condition (4.2) can be given as

$$\begin{bmatrix} (J_1^{(j)})^\mathrm{T} \Gamma^{-1} J_1^{(j)} + \alpha H(\sigma) & (J_1^{(j)})^\mathrm{T} \Gamma^{-1} J_2^{(j)} \\ (J_2^{(j)})^\mathrm{T} \Gamma^{-1} J_1^{(j)} & (J_2^{(j)})^\mathrm{T} \Gamma^{-1} J_2^{(j)} \end{bmatrix} \kappa = (J^{(j)})^\mathrm{T} \Gamma^{-1} y^{(j)}. \qquad (4.5)$$

The next iterate $\kappa^{(j+1)}$ in the sequence produced by the sequential linearizations is defined to be the solution of (4.5). Note that (4.5) is a nonlinear equation for $\kappa = [\sigma^\mathrm{T}, z^\mathrm{T}]^\mathrm{T}$ since $H(\sigma)$ depends on $\sigma$.

**4.2. Lagged diffusivity iteration.** The lagged diffusivity iteration [37] is a fixed point algorithm having its origins in image processing. It is suitable for solving nonlinear equations such as (4.5) and results in a nested iteration inside the sequential linearizations (see Section 5 for a heuristic relief of this computational burden in the final version of our reconstruction algorithm). When estimating $\kappa^{(j+1)}$, an approximating sequence $\kappa_i^{(j+1)}$, $i = 0, 1, \ldots$, is formed by solving $\kappa_{i+1}^{(j+1)}$ from the linear equation obtained from (4.5) by forming $H(\sigma)$ at the conductivity defined by the current iterate $\kappa_i^{(j+1)}$. To be more precise, the lagged diffusivity iteration in its basic form is as follows:

1. Set $\kappa_0^{(j+1)} = \kappa^{(j)}$ and $i = 0$.

2. Solve $\kappa_{i+1}^{(j+1)}$ from the linear equation

$$\begin{bmatrix} (J_1^{(j)})^\mathrm{T}\Gamma^{-1}J_1^{(j)} + \alpha\, H(\sigma_i^{(j+1)}) & (J_1^{(j)})^\mathrm{T}\Gamma^{-1}J_2^{(j)} \\ (J_2^{(j)})^\mathrm{T}\Gamma^{-1}J_1^{(j)} & (J_2^{(j)})^\mathrm{T}\Gamma^{-1}J_2^{(j)} \end{bmatrix} \kappa_{i+1}^{(j+1)} = (J^{(j)})^\mathrm{T}\Gamma^{-1}y^{(j)}. \tag{4.6}$$

3. If the chosen stopping criterion is satisfied, define $\kappa^{(j+1)} = \kappa_{i+1}^{(j+1)}$ and terminate the iteration. Otherwise, set $i \leftarrow i + 1$ and return to step 1.

Let us restructure (4.6) so that the determination of $\sigma = \sigma_{i+1}^{(j+1)}$ (unstable) and that of $z = z_{i+1}^{(j+1)}$ (stable) are separated. We denote

$$B_1 = \Gamma^{-1/2} J_1^{(j)}, \qquad B_2 = \Gamma^{-1/2} J_2^{(j)}, \qquad H = H(\sigma_i^{(j+1)}),$$

where $\Gamma^{-1/2}$ is a Cholesky factor[1] of $\Gamma^{-1}$. (As a single realization of (4.6) is considered, we drop the indices related to the sequential linearizations and the lagged diffusivity iteration in order to retain readability.) Recalling $\kappa = [\sigma^\mathrm{T}, z^\mathrm{T}]^\mathrm{T}$, solving the second line of (4.6) for $z$, and substituting in the first one, we obtain

$$(A^\mathrm{T} A + \alpha H)\sigma = A^\mathrm{T} b. \tag{4.7}$$

Here,

$$A = QB_1 \in \mathbb{R}^{M(M-1)\times N}, \qquad Q = \mathrm{I} - P \in \mathbb{R}^{M(M-1)\times M(M-1)},$$

with I being the identity matrix and $P = B_2(B_2^\mathrm{T} B_2)^{-1} B_2^\mathrm{T}$ the orthogonal projection onto the range $\mathcal{R}(B_2)$ in $\mathbb{R}^{M(M-1)}$. Note that this makes $Q$ the orthogonal projection onto $\mathcal{R}(B_2)^\perp$. Moreover,

$$b = Q\Gamma^{-1/2} y \in \mathbb{R}^{M(M-1)},$$

and

$$z = (B_2^\mathrm{T} B_2)^{-1} B_2^\mathrm{T} (\Gamma^{-1/2} y - B_1 \sigma) \in \mathbb{R}^M. \tag{4.8}$$

The relations $Q^\mathrm{T} = Q$ and $QQ = Q$ are essential for deducing the above formulas.

A couple of remarks are in order. First of all, the main reason for dividing (4.6) into the two parts (4.7)–(4.8) is allowing the use of priorconditioning merely on the conductivity by applying LSQR to a suitably preconditioned version of (4.7). Secondly, the inversion of $B_2^\mathrm{T} B_2 \in \mathbb{R}^{M\times M}$ is relatively inexpensive as the number of electrodes $M$ is reasonably low in practical EIT; according to our experience, $B_2$ is always of full rank, making $B_2^\mathrm{T} B_2$ invertible, which reflects the stability of determining the contact resistances from EIT measurements. In consequence, it is not far-fetched to assume the elements of the coefficient matrix $A^{\mathrm{T}}A + \alpha H$ of (4.7) can be formed explicitly, and after (4.7) is solved, obtaining the whole $\kappa$ via (4.8) is practically free. Finally, notice that solving (4.7) is equivalent to finding the minimizer of the least squares functional

$$\|A\sigma - b\|^2 + \alpha\, \sigma^{\mathrm{T}} H \sigma, \tag{4.9}$$

where $\|\cdot\|$ denotes the Euclidean norm. This can be interpreted as determining the MAP (or conditional mean) estimate corresponding to the linear problem $A\sigma = b$ with a suitable additive Gaussian noise model and a Gaussian mean-free prior with the inverse covariance matrix $H$.

---

[1]In the numerical studies of Section 6, the covariance matrix $\Gamma$ is diagonal, meaning that $\Gamma^{-1/2}$ can be obtained by taking the 'inverse square root' of $\Gamma$. This explains the chosen notation.

**5. Priorconditioning and implementation.** The considerations of Section 4 have three obvious imperfections. First of all, the free parameter $\alpha > 0$ is an extra nuisance. Secondly, we have not yet commented how the problem (4.7) should be solved. Finally, the employment of two nested iterations (sequential linearizations and lagged diffusivity iterations) seems computationally inefficient. In what follows, we explain how these flaws are fixed.

The computationally most expensive part in tackling the system (4.7)–(4.8) is by far solving (4.7) accurately enough for $\sigma \in \mathbb{R}^N$. Here we resort to priorconditioning. We formally (Cholesky) factor $H = L^{\mathrm{T}}L$ and multiply (4.7) from the left by $(L^{-1})^{\mathrm{T}}$, which leads to

$$\big((L^{-1})^{\mathrm{T}} A^{\mathrm{T}} A L^{-1} + \alpha \mathrm{I}\big)(L\sigma) = (L^{-1})^{\mathrm{T}} A^{\mathrm{T}} b. \tag{5.1}$$

Instead of choosing $\alpha$ based on *a priori* knowledge on $\sigma$, we assume that the essential prior information is already encoded in the solution process via the preconditioning by $H = L^{\mathrm{T}}L$ in (5.1) [1]. Hence, we set $\alpha = 0$ and aim at regularizing (5.1) by applying LSQR [29, 28], which is analytically equivalent to the conjugate gradient method for normal equations but works with $A$ instead the more illconditioned $A^{\mathrm{T}}A$, combined with a suitable *early stopping rule* [17, 19]. We assume the knowledge of an approximate constant background conductivity level $\sigma^{(0)}$ and rewrite (5.1), with $\alpha = 0$, in the form

$$(L^{-1})^{\mathrm{T}} A^{\mathrm{T}} A L^{-1} \tilde{\sigma} = (L^{-1})^{\mathrm{T}} A^{\mathrm{T}} \tilde{b}, \qquad \sigma = L^{-1}\tilde{\sigma} + \sigma^{(0)}, \tag{5.2}$$

where $\tilde{b} = b - A\sigma^{(0)}$. The aim of the translation by $\sigma^{(0)}$ is to allow implementing a simple homogeneous Dirichlet boundary condition for (4.4) to guarantee the invertibility of $H$; cf. Algorithm 1 below. The equation (5.2) can be solved efficiently using the modified LSQR algorithm from [1]. In particular, the factorization $H = L^{\mathrm{T}}L$ is introduced *merely for notational convenience*: The LSQR algorithm from [1] does not form such. However, one needs to be able to apply the inverse of $H$ on a given vector.

The reason for employing LSQR to the preconditioned system (5.2) instead of (4.7) is the form of the associated Krylov subspaces. For the initial guess $\tilde{\sigma} = 0$, the approximation of $\sigma$ produced by LSQR after $m$ iterations belongs to [1]

$$\mathcal{K}_m = \mathrm{span}\{H^{-1}A^{\mathrm{T}}\tilde{b}, (H^{-1}A^{\mathrm{T}}A)H^{-1}A^{\mathrm{T}}\tilde{b}, \dots, (H^{-1}A^{\mathrm{T}}A)^{m-1}H^{-1}A^{\mathrm{T}}\tilde{b}\}$$

shifted by $\sigma^{(0)}$. Thus, the approximate solution is of the form

$$\sigma = H^{-1}a + \sigma^{(0)} \tag{5.3}$$

for some $a \in \mathbb{R}^N$ independently of the number of iterations. An operation with $H \in \mathbb{R}^{N \times N}$ can be considered to measure the compatibility of a given vector with the prior information: The smaller the squared $H$-norm $\sigma^{\mathrm{T}} H \sigma$ is, the better is the correspondence between $\sigma$ and the (linearized) prior (cf. (4.9)). Inverting this logic, vectors of the form (5.3) should contain features that are in accordance with the prior. We refer to [1, 8, 25] for more information on this phenomenon.

Since (5.2) lacks a penalty term, an early stopping rule must be employed. Our choice is the Morozov discrepancy principle [17, 19]: The LSQR iteration is stopped when the squared residual for the least squares problem corresponding to (4.6) equals $\mathbb{E}(\|\Gamma^{-1/2}\eta\|^2) = M(M-1)$, that is,

$$\|B\kappa - \Gamma^{-1/2}y\|^2 = M(M-1), \tag{5.4}$$

where $B = [B_1, B_2]$ and we have once again dropped the indices for the sake of readability. The condition (5.4) can easily be monitored when applying LSQR to (5.1) because it holds that

$$\begin{aligned}\|B\kappa - \Gamma^{-1/2}y\|^2 &= \|B_1\sigma + B_2 z - \Gamma^{-1/2}y\|^2 \\ &= \|B_1\sigma + B_2(B_2^{\mathrm{T}} B_2)^{-1} B_2^{\mathrm{T}}(\Gamma^{-1/2}y - B_1\sigma) - \Gamma^{-1/2}y\|^2 \\ &= \|QB_1\sigma - Q\Gamma^{-1/2}y\|^2 = \|A\sigma - b\|^2,\end{aligned}$$

where the second step follows by substituting (4.8). As a consequence, the Morozov discrepancy principle can be applied directly to LSQR, that is, the iteration is stopped when $\|A\sigma - b\|$ equals $\sqrt{M(M-1)}$. This observation has an intuitive explanation: The orthogonal projection of the data $\Gamma^{-1/2}y$ onto $\mathcal{R}(B_2)$ is completely explained by the choice (4.8) for $z$, and thus the discrepancy in (5.4) originates solely from the projection of the considered equation onto the orthogonal complement $\mathcal{R}(B_2)^{\perp}$.

As for the two nested iterations — or three, if LSQR is counted in —, our experience suggests that it is enough to take just a single lagged diffusivity step for each linearization of the CEM forward operator, that is, for each realization of (4.5). A higher number of iterations does not seem to have any significant effect on the final reconstructions, leading only to longer computation times. In other words, the final version of our algorithm changes the basis point for the linearization of $(\sigma, z) \mapsto \mathcal{U}(\sigma, z)$ after each lagged diffusivity update. Moreover, take note that there is no sensible reason for not updating the linearization after each lagged diffusivity step if the reconstruction algorithm is implemented in a matrix-free manner (, which is not the case in our numerical experiments but anyway a possibility [1]).

Including a global Morozov-type stopping criterion for the whole reconstruction process, our algorithm is altogether as follows.

ALGORITHM 1. *Pick $T > 0$ and your favorite $r$, e.g., from (3.5). Choose $(\sigma^{(0)}, z^{(0)}) \in \mathbb{R}^N \times \mathbb{R}^M$ to be the minimizer of*

$$E(\sigma, z) := \|\Gamma^{-1/2}(\mathcal{V} - \mathcal{U}(\sigma, z))\| \tag{5.5}$$

*over the pairs of positive homogeneous conductivities and contact resistance vectors (two free parameters). Set $j = 0$, $\epsilon = \sqrt{M(M-1)}$ and $\mathcal{U} = \mathcal{U}(\sigma^{(0)}, z^{(0)})$.*

1. *Evaluate the Jacobian matrices of the map $(\sigma, z) \mapsto \mathcal{U}(\sigma, z)$ with respect to $\sigma$ and $z$ at $(\sigma^{(j)}, z^{(j)})$. Denote these entities as $J_1$ and $J_2$, respectively.*

2. *Set $y = \mathcal{V} - \mathcal{U} + J_1 \sigma^{(j)} + J_2 z^{(j)}$, $B_1 = \Gamma^{-1/2} J_1$, $B_2 = \Gamma^{-1/2} J_2$, and let $Q$ be the orthogonal projection onto $\mathcal{R}(B_2)^{\perp}$. Define $A = QB_1$ and $\tilde{b} = Q(\Gamma^{-1/2}y - B_1\sigma^{(0)})$.*

3. Build the preconditioner $H = H(\sigma^{(j)})$ according to (4.3), including a homogeneous Dirichlet boundary condition for (4.4) on the electrodes $E$ to guarantee invertibility.

4. Apply the LSQR algorithm from [1] to
$$(L^{-1})^{\mathrm{T}} A^{\mathrm{T}} A L^{-1} \tilde{\sigma} = (L^{-1})^{\mathrm{T}} A^{\mathrm{T}} \tilde{b}, \qquad \sigma = L^{-1}\tilde{\sigma} + \sigma^{(0)},$$
starting from $\tilde{\sigma} = 0$, which results in a solution sequence $\{\sigma_k^{(j)}\}_{k\geq 0}$. Stop the iteration when $\|A\sigma_k^{(j)} - b\| \leq \epsilon$. Denote the corresponding solution as $\sigma^{(j+1)}$. (Here, $L$ satisfies $H = L^{\mathrm{T}} L$, but it need not be formed explicitly [1].)

5. Set $z^{(j+1)} = (B_2^{\mathrm{T}} B_2)^{-1} B_2^{\mathrm{T}} (\Gamma^{-1/2} y - B_1 \sigma^{(j+1)})$.

6. Compute $\mathcal{U} = \mathcal{U}(\sigma^{(j+1)}, z^{(j+1)})$. If $E(\sigma^{(j+1)}, z^{(j+1)}) \leq \epsilon$, dub $(\sigma^{(j+1)}, z^{(j+1)})$ the reconstruction and exit. Otherwise, set $j \leftarrow j+1$ and return to step 1.

Assuming the noise covariance matrix $\Gamma$ is known, Algorithm 1 has only one free parameter, namely $T$. In all numerical experiments of Section 6, we use $T = 10^{-6}$ for total variation and $T = 5 \times 10^{-3}$ for Perona–Malik but these choices do not have a significant effect on the reconstructions.

Notice that the initialization step of choosing $(\sigma^{(0)}, z^{(0)})$ consists of determining just two positive real parameters: the constant conductivity level and the universal contact resistance that together best explain the data. Subsequently, the initial constant conductivity level $\sigma^{(0)}$ is used implicitly as a Dirichlet boundary value for all succeeding $\sigma^{(j)}$, $j = 1, 2, \ldots$, on $E \subset \partial\Omega$. Indeed, since $H$ is accompanied by a homogeneous Dirichlet boundary condition on $E$ in step 3 of Algorithm 1, it holds that $\sigma|_E = \sigma^{(0)}|_E$ due to (5.3). The reason for selecting the electrodes as the boundary part with the homogeneous Dirichlet boundary condition for (4.4) is the ability of the contact resistances to compensate for slightly false conductivity underneath the electrodes. Be that as it may, the choice of the boundary region with the Dirichlet boundary condition for (4.4) does not seem to have a huge effect on the reconstructions.

To conclude this section, let us briefly comment on the computational cost of applying the inverse of $H \in \mathbb{R}^{N \times N}$ on a given vector in each step of the LSQR iterations. Although the number of FEM nodes $N$ is relatively high in our three-dimensional numerical experiments on dense grids, the sparse structure of $H$ makes it possible to just use the 'backslash' operation in MATLAB without running into trouble regarding the computation time. (On the other hand, note that all Jacobian matrices appearing in Algorithm 1 are full.) However, if the number of degrees of freedom $N$ were even higher, a reasonable approach would be to employ a suitable multigrid solver [4] to (approximately) multiply with $H^{-1}$. In [1] it was demonstrated that just a couple of V-cycles of a suitable multigrid method is enough for the purpose of priorconditioning in the case of a related linear inverse problem; according to our experience, the same applies to the numerical tests of Section 6 as well.

**6. Numerical experiments.** We apply Algorithm 1 to three test cases. The first two consider simulated data and generic geometrical settings with several electrodes, while the third one deals with experimental data from a thorax-shaped water tank with sixteen electrodes. In the first two experiments, the electrode potential measurements are simulated via solving (2.1) by FEM with piecewise linear basis functions (cf. [34, 35]) on a considerably denser finite element mesh than the one em-

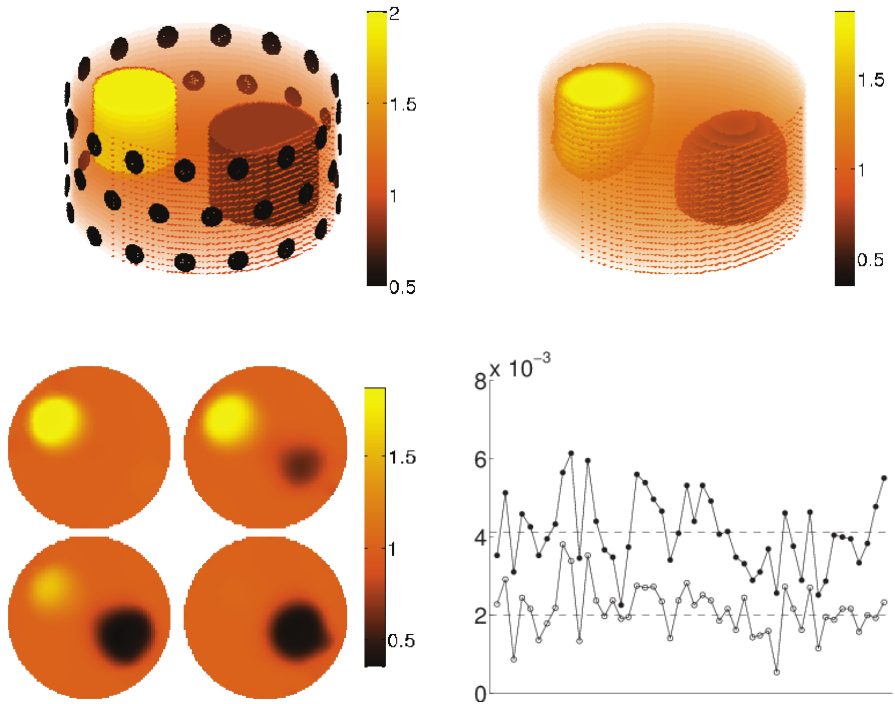

Fig. 6.1: Results for Case 1 with Perona–Malik. Top left: The target conductivity and the electrodes. Top right: The conductivity reconstruction; conductivities in the interval $[0.8, 1.2]$ are transparent. Bottom left: The conductivity reconstruction; four cross-sections at heights $0.9, 0.6, 0.4$ and $0.1$. Bottom right: The reconstructed contact resistances (filled) compared with the exact ones. The dashed lines indicate the respective mean values.

ployed in Algorithm 1. Subsequently, the actual data vector $\mathcal{V}$ is formed according to (3.2) with

$$\Gamma = \left(\frac{\gamma}{M(M-1)} \sum_{j=1}^{M(M-1)} |\mathcal{U}_j(\sigma_{\text{xct}}, z_{\text{xct}})|\right)^2 \mathrm{I} \in \mathbb{R}^{M(M-1)\times M(M-1)}, \qquad (6.1)$$

where $\mathcal{U}(\sigma_{\text{xct}}, z_{\text{xct}}) \in \mathbb{R}^{M(M-1)}$ is the exact (simulated) measurement vector and $\gamma > 0$ is a free parameter. In other words, the noise covariance matrix is assumed to be diagonal and homogeneous; the components of the additive measurement noise vector are independent but have the same standard deviation, which is $100\gamma$ percent of the mean *absolute* electrode potential measurement. In all three tests, the applied $M-1$ linearly independent current patterns are $I^m = e_1 - e_{m+1}$, $m = 1, \ldots, M-1$, where $e_m$ denotes the $m$th Cartesian basis vector, meaning that one feeding electrode is fixed and the current exits in turns through the remaining ones.

**Case 1: Cylindrical body (simulated data).** The set-up of the first numerical test is illustrated in the top left image of Figure 6.1. The examined object is a cylinder $\Omega = D(0,1) \times (0,1) \subset \mathbb{R}^3$, with $D(x,r) \subset \mathbb{R}^2$ denoting an open disk of radius $r > 0$ centered at $x \in \mathbb{R}^2$. There are $M = 48$ identical, disk-shaped electrodes attached to the side of $\Omega$. The internal conductivity of $\Omega$ is $\sigma \equiv 1$ apart from two embedded disk-based cylindrical inclusions of height 0.6 and radii 0.4 and 0.3, respectively, with axes parallel to that of $\Omega$. One of the inclusions stands at the bottom of $\Omega$ and has conductivity $\sigma = 0.5$, the other 'hangs' from the top of $\Omega$ and is of conductivity $\sigma = 2$. The target contact resistances are of the form

$$z_m = \tilde{z} + v_m, \qquad m = 1, \ldots, M, \tag{6.2}$$

where $\tilde{z} = 2 \times 10^{-3}$ and $v_m$, $m = 1, \ldots, M$, are independent realizations of a Gaussian random variable with vanishing mean and standard deviation $5 \times 10^{-4}$ (cf. Figure 6.1). The parameter controlling the amount of measurement noise in (6.1) is chosen to be $\gamma = 4 \times 10^{-3}$, which loosely speaking corresponds to 0.4% of noise.

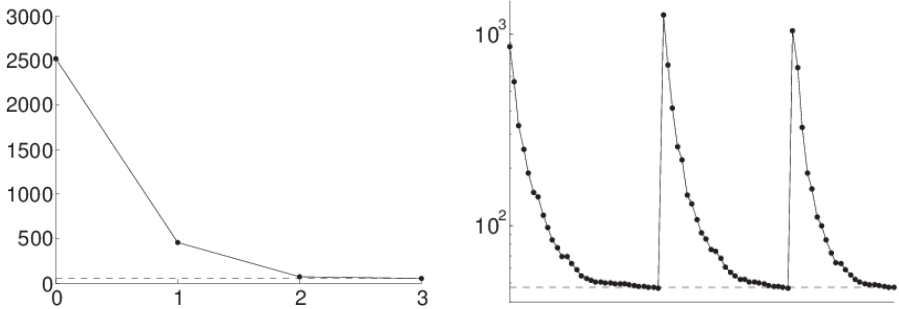

Fig. 6.2: Results for Case 1 with Perona-Malik. Left: The evolution of the residual $E(\sigma^{(j)}, z^{(j)})$ as a function of $j$. Right: The 'linearized' residual $\|A\sigma_k^{(j)} - b\|$ on a logarithmic scale over both the sequential linearizations $j$ and the LSQR iterations $k$. The dashed lines indicate the noise level $\epsilon$.

We start with the Perona–Malik prior, i.e., the second option for $r$ in (3.5). The conductivity reconstruction produced by Algorithm 1 is visualized by two different techniques in the top right and bottom left images of Figure 6.1. The corresponding contact resistances are compared with the exact ones in the bottom right image. Although the reconstructed inclusions do not have precisely the correct shapes, they are well localized, lie close to the correct positions and reproduce the size and conductivity levels of the target inhomogeneities relatively accurately. Moreover, the background in the reconstruction of Figure 6.1 is approximately constant. Hence, Algorithm 1 seems to promote the desired qualities. Notice that the reconstructed contact resistances are systematically too high, which is a feature of Algorithm 1 that recurs throughout our numerical studies. A possible reason is the decoupling of the conductivity and contact resistances in (4.7)–(4.8), which may allow the contact resistances to explain part of the internal resistivity of the object. Be that as it may, this mild flaw does not seem to severely impair the conductivity reconstruction, i.e., the main objective of EIT.

The left-hand image of Figure 6.2 presents the residual $E(\sigma^{(j)}, z^{(j)})$, defined in (5.5), as a function of $j$, i.e., the number of outer iterations in Algorithm 1, whereas the right-hand image shows the behavior of the 'linearized residual' $\|A\sigma_k^{(j)} - b\|$ over both the sequential linearizations $j$ and the LSQR iteration $k$. It is worth noticing that as $j$ increases the number of LSQR steps needed for reaching the noise level $\epsilon = \sqrt{M(M-1)}$ decreases. This behavior, which is in line with the observations in [1] and consistent throughout our numerical studies, increases the computational

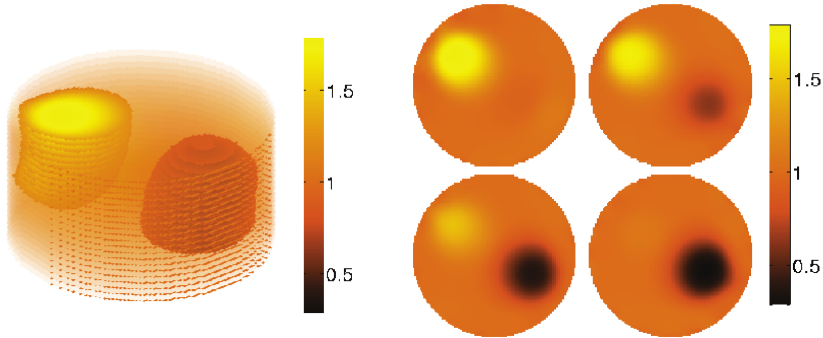

Fig. 6.3: Results for Case 1 with total variation. Left: The conductivity reconstruction; the conductivities in the interval $[0.8, 1.2]$ are transparent. Right: The conductivity reconstruction; four cross-sections at heights $0.9, 0.6, 0.4$ and $0.1$.

efficiency of the algorithm.

We complete the discussion on Case 1 by repeating the process with the total variation prior, i.e., the first option in (3.5). The resulting conductivity reconstruction is visualized in Figure 6.3. Although Algorithm 1 functions relatively well also in this case, according to visual inspection, the reconstruction in Figure 6.1 is somewhat better than the one in Figure 6.3. As a consequence, we use Perona–Malik in the remaining tests; however, the difference between the reconstructions corresponding to the two options in (3.5) is not significant in Cases 2 and 3 either. It is also worth mentioning that by using the so-called TV$^q$ prior, which is analogous to total variation but approximates $|t|^q$, $0 < q < 1$, instead of $|t|$, we get almost as good reconstructions as with Perona–Malik.

The reconstructions presented in Figures 6.1 and 6.3 were formed on a finite element mesh with 80300 tetrahedrons and 17031 nodes, with appropriate refinements at the edges of the electrodes. For Perona–Malik, the running time of the algorithm was 177 seconds with a MATLAB-based implementation on a desktop with 8GB RAM and an Intel Xeon E31230 CPU having clock speed 3.20 GHz. For total variation, the computation time was 178 seconds.

**Case 2: Ball (simulated data).** The top left image of Figure 6.4 shows the configuration of the second test case. This time the object of interest $\Omega$ is a ball of radius 5 and the measurements are performed with $M = 37$ identical disk-shaped electrodes. The body $\Omega$ has a constant background conductivity level $\sigma = 1$ with two embedded inclusions of conductivities 0.5 and 2, respectively. One of the inhomogeneities is the intersection of $\Omega$ and a ball of radius 2; in particular, it touches the boundary $\partial\Omega$. The other is L-shaped and, in particular, nonconvex. The contact resistances are of the form (6.2) with $\tilde{z} = 0.01$ and the standard deviation 0.002 for $v_m$, $m = 1, \ldots, M$. The measurement noise level in (6.1) is chosen again as $\gamma = 4 \times 10^{-3}$.

The conductivity reconstruction provided by Algorithm 1 with the Perona–Malik prior is visualized in the top right and bottom left images of Figure 6.4. The conclusions are the same as in Case 1: The reconstructed inclusions are well localized, have approximately the correct size and conductivity level, and are embedded in a homogeneous background. In particular, even the inclusion attached to the boundary is reconstructed reasonably well. The second horizontal cross-section shows the reconstruction at the level of the nonconvex inhomogeneity: Although the shape of the inclusion is not reproduced accurately, the nonconvexity is retained in the reconstruction. In this test, the 'nonlinear' residual $E(\sigma^{(j)}, z^{(j)})$ did not reach the noise level $\epsilon$ before the reconstruction became unstable. Hence, the global Morozov parameter was increased to $1.5\epsilon$, while the stopping criterion for the individual LSQR iterations remained as previously. The bottom right image of Figure 6.4 shows the behavior of the 'linearized residual' $\|A\sigma_k^{(j)} - b\|$ over both the sequential linearizations $j$ and the LSQR iteration $k$.

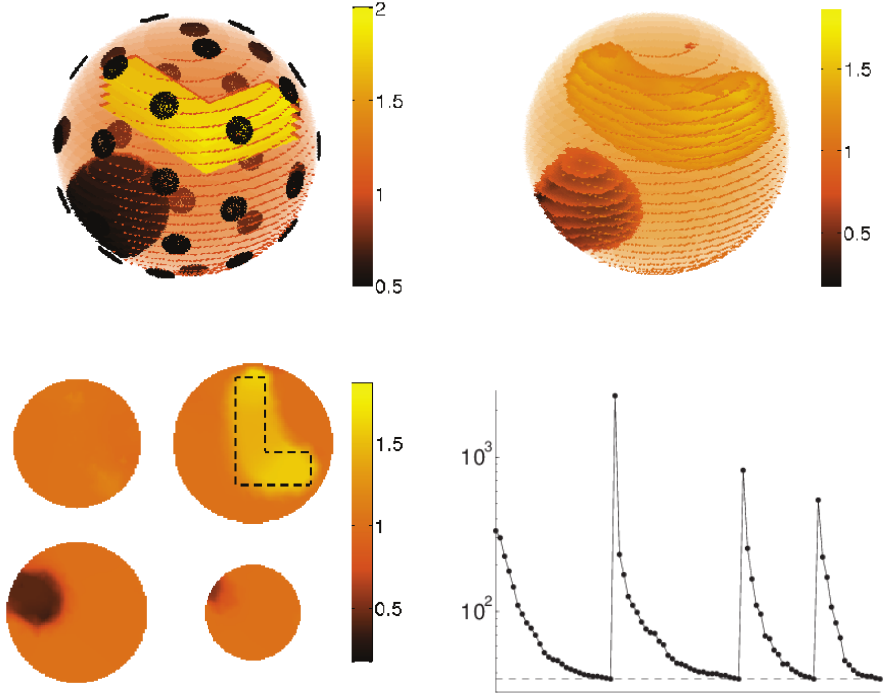

Fig. 6.4: Results for Case 2 with Perona–Malik. Top left: The target conductivity and the electrodes. Top right: The conductivity reconstruction; conductivities in the interval $[0.8, 1.2]$ are transparent. Bottom left: The conductivity reconstruction; four cross-sections at heights $3, 1, -2.5$ and $-4$. The dashed curve at height $1$ indicates the exact inclusion boundary. Bottom right: The 'linearized' residual $\|A\sigma_k^{(j)} - b\|$ on a logarithmic scale over both the sequential linearizations $j$ and the LSQR iterations $k$. The dashed line indicates the noise level $\epsilon$.

The reconstruction presented in Figure 6.4 was computed on a finite element mesh with 45896 tetrahedrons and 11608 nodes. The total number of LSQR iterations was $27 + 29 + 17 + 15 = 88$ and the running time of Algorithm 1 was 86 seconds with the same implementation and hardware as in Case 1.

**Case 3: Water tank.** Our third and final test case deals with experimental data from a thorax-shaped water tank. The tank has circumference $106\,\mathrm{cm}$ and there are $M = 16$ rectangular electrodes of width $2\,\mathrm{cm}$ and height $5\,\mathrm{cm}$ attached to its interior lateral surface. We consider the two inclusion configurations shown in the top left images of Figures 6.5 and 6.6: The former considers one insulating cylindrical inclusion made of plastic whereas the latter corresponds to one plastic and one highly conductive metallic inclusion cylinder. The tank is filled with Finnish tap water up to the top of the electrodes. The measurements were performed with low-frequency alternating current (frequency 1 kHz, amplitude 1 mA) using the *Kuopio impedance tomography* (KIT4) device [27]. The phase information of the voltage measurements was ignored, and the electrode currents and potentials were interpreted as real vectors, with the aim being the reconstruction of an approximate real admittivity, i.e., conductivity.

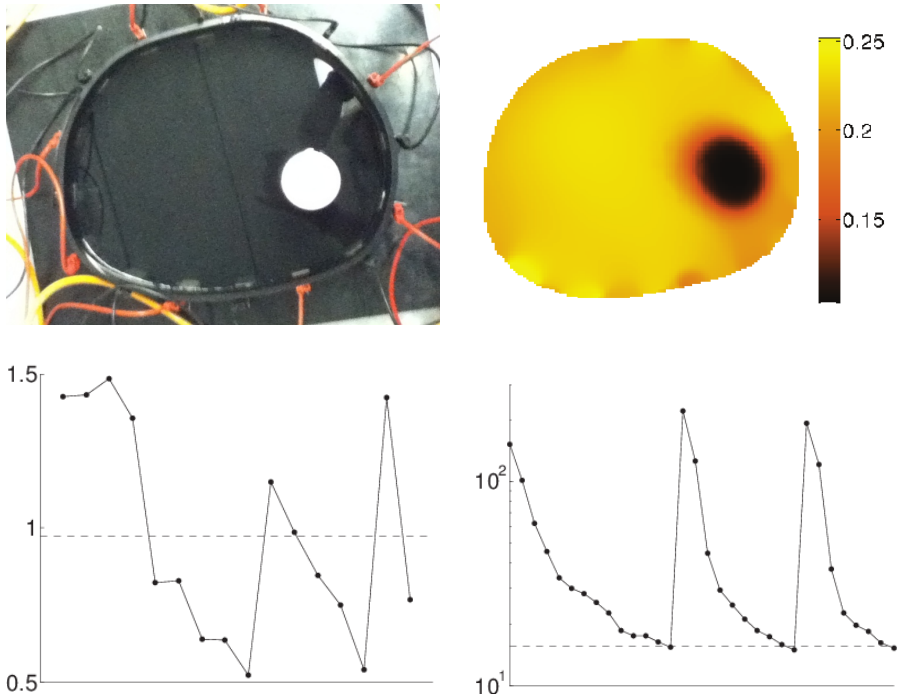

Fig. 6.5: Results for Case 3 with Perona–Malik. Top left: The target. Top right: The cross-section of the conductivity reconstruction (mS/cm) at height 2.5 cm. Bottom left: The reconstructed contact resistances (k$\Omega$ cm$^2$). The dashed line indicates the mean value. Bottom right: The 'linearized' residual $\|A\sigma_k^{(j)} - b\|$ on a logarithmic scale over both the sequential linearizations $j$ and the LSQR iterations $k$. The dashed line indicates the noise level $\epsilon$.

Although our vertically homogeneous measurement setting could be modelled by a two-dimensional version of the CEM, our computations are performed on a three-dimensional FEM mesh with 136286 tetrahedrons and 27276 nodes. This time around we choose the form

$$\Gamma = \left(\gamma \max_{i,j} |\mathcal{V}_i - \mathcal{V}_j|\right)^2 \mathrm{I} \in \mathbb{R}^{M(M-1) \times M(M-1)} \tag{6.3}$$

for the noise covariance matrix, that is, the standard deviation of the additive Gaussian measurement noise on each electrode is assumed to be proportional to the difference of the smallest and largest electrode potential measurement. By trial and error, the free parameter in (6.3) was set as $\gamma = 6 \times 10^{-4}$; such a choice results in good reconstructions and a reasonable speed of convergence for the combination of

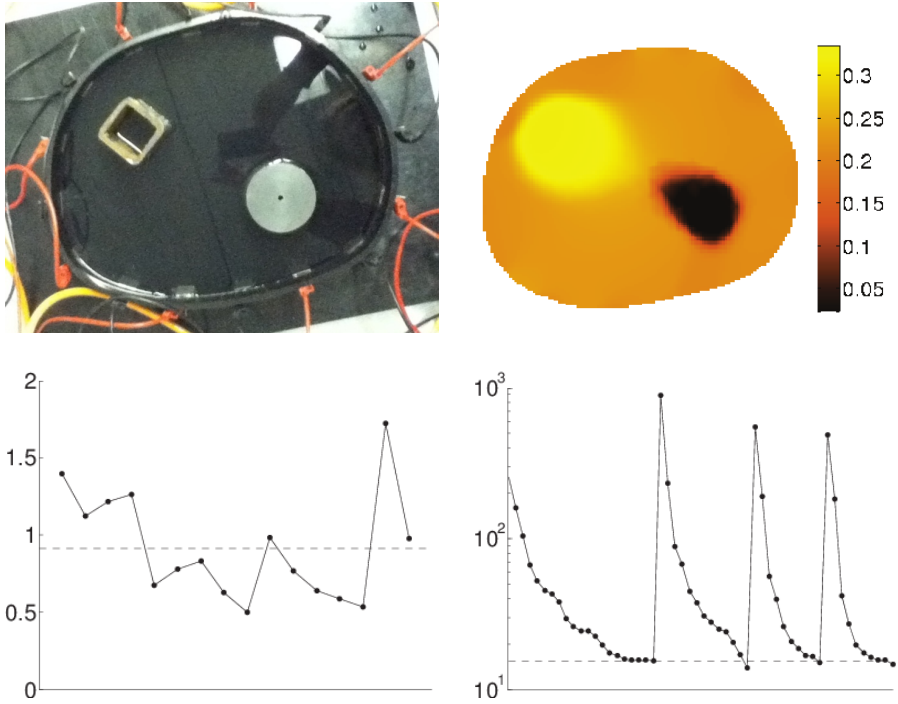

Fig. 6.6: Results for Case 3 with Perona–Malik. Top left: The target. Top right: The cross-section of the conductivity reconstruction (mS/cm) at height $2.5\,\text{cm}$. Bottom left: The reconstructed contact resistances ($\text{k}\Omega\,\text{cm}^2$). The dashed line indicates the mean value. Bottom right: The 'linearized' residual $\|A\sigma_k^{(j)}-b\|$ on a logarithmic scale over both the sequential linearizations $j$ and the LSQR iterations $k$. The dashed line indicates the noise level $\epsilon$.

Algorithm 1 and the KIT4 device.

The reconstructions produced by Algorithm 1 with Perona–Malik prior for the two target conductivities are illustrated in Figures 6.5 and 6.6, which are organized in the same way. In each case, the top right image shows the horizontal cross-section of the reconstruction at height $2.5\,\text{cm}$. The bottom left image depicts the corresponding contact resistances, which are considerably higher than the values encountered in previous water tank experiments (cf., e.g., [18]). The bottom right image shows the behavior of the 'linearized residual' $\|A\sigma_k^{(j)}-b\|$ over both the sequential linearizations $j$ and the LSQR iteration $k$. Like for simulated CEM measurements in Cases 1 and 2, Algorithm 1 produces good reconstructions of the target conductivity also from experimental data. The shapes of the reconstructed inclusions are not perfect, but their positions are reproduced accurately. Moreover, the inhomogeneities are well localized and the background conductivity does not fluctuate considerably, apart from some regions close to the object boundary. This undesirable effect is probably due to the unavoidable mismodelling of the water tank shape, which could possibly be compensated by including the fine-tuning of the geometric information as a part of the reconstruction algorithm (cf. [13, 14]). Notice that the conductivity of the metallic cylinder in Figure 6.6 is estimated only as $0.33\,\text{mS/cm}$, which is not necessarily only a demerit of the algorithm but can also be caused by contact resistance on the surface of the inclusion in question; see, e.g., [23] and the references therein.

The total number of LSQR iterations and the running time of Algorithm 1 were $14 + 10 + 8 = 32$ and 146 seconds, respectively, for the setting of Figure 6.5 and $21 + 13 + 10 + 10 = 54$ and 180 seconds, respectively, for that of Figure 6.6. The implementation and the hardware were the same as in Cases 1 and 2.

**7. Conclusion.** We have introduced a numerical algorithm for EIT capable of producing useful three-dimensional reconstructions of inclusion-type target conductivities on dense FEM meshes in computation times of only a few minutes. The method is based on the use of a nonlinear edge-enhancing penalty term and the minimization of the corresponding Tikhonov functional by efficiently solving an approximate sequence of linearized problems with the help of priorconditioning and LSQR. The algorithm was successfully tested both with simulated CEM measurements and with experimental water tank data.

**Acknowledgments.** We would like to thank Professor Jari Kaipio's research group at the University of Eastern Finland (Kuopio) for letting us use their three-dimensional FEM solver for the CEM as well as for granting us access to their EIT devices. We are particularly grateful to Dr. Aku Seppänen.

REFERENCES

[1] ARRIDGE, S., BETCKE, M., AND HARHANEN, L. Iterated preconditioned LSQR method for inverse problems on unstructured grids. *Inverse Problems* (2014). Accepted.
[2] BORCEA, L. Electrical impedance tomography. *Inverse problems 18* (2002), R99–R136.
[3] BORSIC, A., GRAHAM, B. M., ADLER, A., AND LIONHEART, W. R. B. In vivo impedance imaging with total variation regularization. *IEEE Trans. Med. Imaging 29* (2010), 44–54.
[4] BRIGGS, W. L., VAN EMDEN, H., AND MCCORMICK, S. F. *A multigrid tutorial*. SIAM, 2000.
[5] CALVETTI, D. Preconditioned iterative methods for linear discrete ill-posed problems from a bayesian inversion perspective. *J. Comput. Appl. Math. 198* (2007), 378395.
[6] CALVETTI, D., MCGIVNEY, D., AND SOMERSALO, E. Left and right preconditioning for electrical impedance tomography with structural information. *Inverse Problems 28* (2012), 055015.
[7] CALVETTI, D., AND SOMERSALO, E. Priorconditioners for linear systems. *Inverse problems 21* (2005), 1397–1418.
[8] CALVETTI, D., AND SOMERSALO, E. *Introduction to Bayesian scientific computing: Ten lectures on subjective computing*. Springer, 2007.
[9] CALVETTI, D., AND SOMERSALO, E. Hypermodels in the bayesian imaging framework. *Inverse Problems 24* (2008), 034013.
[10] CHENEY, M., ISAACSON, D., AND NEWELL, J. Electrical impedance tomography. *SIAM Rev. 41* (1999), 85–101.
[11] CHENG, K.-S., ISAACSON, D., NEWELL, J. S., AND GISSER, D. G. Electrode models for electric current computed tomography. *IEEE Trans. Biomed. Eng. 36* (1989), 918–924.
[12] CHUNG, E. T., CHAN, T. F., AND TAI, X.-C. Electrical impedance tomography using level set representation and total variational regularization. *J. Comput. Phys. 205* (2005), 357372.
[13] DARDÉ, J., HYVÖNEN, N., SEPPÄNEN, A., AND STABOULIS, S. Simultaneous reconstruction of outer boundary shape and admittance distribution in electrical impedance tomography. *SIAM J. Imaging Sci. 6* (2013), 176–198.
[14] DARDÉ, J., HYVÖNEN, N., SEPPÄNEN, A., AND STABOULIS, S. Simultaneous recovery of admittivity and body shape in electrical impedance tomography: An experimental evaluation. *Inverse Problems 29* (2013), 085004.
[15] DARDÉ, J., AND STABOULIS, S. Electrode modelling: The effect of contact impedance. arXiv:1312.4202.
[16] ELDÉN, L. A weighted pseudoinverse, generalized singular values, and constrained least squares problem. *BIT 22* (1982), 487–502.
[17] ENGL, H., HANKE, M., AND NEUBAUER, A. *Regularization of inverse problems*. Kluwer, 1996.
[18] GEHRE, M., KLUTH, T., LIPPONEN, A., JIN, B., SEPPÄNEN, A., KAIPIO, J. P., AND MAASS, P. Sparsity reconstruction in electrical impedance tomography: An experimental evaluation. *J. Comp. Appl. Math. 236* (2012), 2126–2136.

[19] HANKE, M. *Conjugate gradient type methods for ill-posed problems*. Chapman & Hall, 1995.
[20] HANSEN, P. C. *Rank-Deficient and Discrete Ill-Posed Problems*. SIAM, 1998.
[21] HARRACH, B., AND SEO, J. K. Exact shape-reconstruction by one-step linearization in electrical impedance tomography. *SIAM J. Math. Anal. 42* (2010), 1505–1518.
[22] HILGERS, J. W. On the equivalence of regularization and certain reproducing kernel hilbert space approaches for solving first kind problems. *SIAM J. Numer. Anal. 13* (1976), 172–184.
[23] HYVÖNEN, N., SEPPÄNEN, A., AND KARHUNEN, K. Fréchet derivative with respect to the shape of an internal electrode in electrical impedance tomography. *SIAM J. Appl. Math. 70* (2010), 1878–1898.
[24] KAIPIO, J. P., KOLEHMAINEN, V., SOMERSALO, E., AND VAUHKONEN, M. Statistical inversion and Monte Carlo sampling methods in electrical impedance tomography. *Inverse Problems 16* (2000), 1487–1522.
[25] KILMER, M., HANSEN, P., AND ESPANOL, M. A projection based approach to general form tikhonov regularization. *SIAM J. Sci. Comput. 29* (2007), 315–330.
[26] KOLEHMAINEN, V., SOMERSALO, E., VAUHKONEN, P. J., VAUHKONEN, M., AND KAIPIO, J. P. A Bayesian approach and total variation priors in 3D electrical impedance tomography. *in Proceedings of the 20th Annual International Conference of the IEEE Engineering in Medicine and Biology Society 20* (1998), 1028–1031.
[27] KOURUNEN, J., SAVOLAINEN, T., LEHIKOINEN, A., VAUHKONEN, M., AND HEIKKINEN, L. M. Suitability of a PXI platform for an electrical impedance tomography system. *Meas. Sci. Technol. 20* (2009), 015503.
[28] PAIGE, C. C., AND SAUNDERS, M. A. Algorithm 583: LSQR: Sparse linear equations and least squares problems. *ACM Trans. Math. Softw. 8* (1982), 195–209.
[29] PAIGE, C. C., AND SAUNDERS, M. A. LSQR: An algorithm for sparse linear equations and sparse least squares. *ACM Trans. Math. Softw. 8* (1982), 43–71.
[30] PERONA, P., AND MALIK, J. Scale-space and edge detection using anisotropic diffusion. *IEEE T. Pattern Anal. 12* (1990), 629639.
[31] RUDIN, L. I., OSHER, S., AND FATEMI, E. Nonlinear total variation based noise removal algorithms. *Physica D 60* (1992), 259–268.
[32] SOMERSALO, E., CHENEY, M., AND ISAACSON, D. Existence and uniqueness for electrode models for electric current computed tomography. *SIAM J. Appl. Math. 52* (1992), 1023–1040.
[33] UHLMANN, G. Electrical impedance tomography and Calderón's problem. *Inverse Problems 25* (2009), 123011.
[34] VAUHKONEN, M. *Electrical impedance tomography with prior information*, vol. 62. Kuopio University Publications C (Dissertation), 1997.
[35] VAUHKONEN, P. J., VAUHKONEN, M., SAVOLAINEN, T., AND KAIPIO, J. P. Three-dimensional electrical impedance tomography based on the complete electrode model. *IEEE Trans. Biomed. Eng. 46* (1999), 11501160.
[36] VILHUNEN, T., KAIPIO, J. P., VAUHKONEN, P. J., SAVOLAINEN, T., AND VAUHKONEN, M. Simultaneous reconstruction of electrode contact impedances and internal electrical properties: I. Theory. *Meas. Sci. Technol. 13* (2002), 18481854.
[37] VOGEL, C. R., AND OMAN, M. E. Iterative methods for total variation denoising. *SIAM J. Sci. Comput. 17* (1996), 227238.